\documentclass[12pt]{amsart}
\usepackage{latexsym,amsmath, amscd, amssymb, amsthm, euscript,mathrsfs, multicol}
\newtheorem{thm}[subsection]{Theorem}
\newtheorem{cor}[subsection]{Corollary}
\newtheorem{lem}[subsection]{Lemma}

\newtheorem{prop}[subsection]{Proposition}

\theoremstyle{definition}

\theoremstyle{definition}

\theoremstyle{definition}

\newtheorem{rem}[subsection]{Remark}

\numberwithin{equation}{subsection}

\newtheorem{pg}[subsection]{}

\newcommand{\mc}{\mathcal }

\newcommand{\Z}{\mathbb{Z}}

\newcommand{\Sp}{\text{\rm Spec}}
\newcommand{\mymargin}[1]{}%{\marginpar{\tiny{#1}}}

\newcommand{\mls}{\mathscr}
\begin{document}

\title{The Picard group of $\mls M_{1,1}$}

\author{William Fulton and Martin Olsson}

\begin{abstract}We compute the Picard  group of the moduli stack of elliptic curves and its canonical compactification over general base schemes
\end{abstract}

\maketitle

\section{Introduction}

Let $\mls  M _{1,1}$ denote the moduli stack (over $\Z$) classifying elliptic curves, and for a scheme $S$ let $\mls M_{1,1,S}$ denote the fiber product $S\times _{\Sp (\Z)}\mls M_{1,1}$.  In his  1965 paper \cite{Mumford}, Mumford computed the Picard group $\text{Pic}(\mls M_{1,1,S})$ when $S$ is the spectrum of a field of characteristic not $2$ or $3$ and found it to be cyclic of order $12$.  Our aim in this  paper is to compute the Picard group $\text{Pic}(\mls M_{1,1,S})$ for more general base schemes $S$, as well as to compute the Picard group $\text{Pic}(\overline {\mls M}_{1,1,S})$ for the standard compactification $\overline {\mls M}_{1,1}$ of $\mls M_{1,1}$.

Recall that on $\mls M_{1,1}$ there is the Hodge bundle $\lambda $.  For any morphism $t:T\rightarrow \mls M_{1,1}$ corresponding to an elliptic curve $f:E\rightarrow T$ the pullback $t^*\lambda $ is the line bundle $f_*\Omega ^1_{E/T}$. Equivalently, if $f:\mls E\rightarrow \mls M_{1,1}$ denotes the universal elliptic curve then $\lambda = f_*\Omega ^1_{\mls E/\mls M_{1,1}}$.  This bundle extends canonically to $\overline {\mls M}_{1,1}$.   Namely, let $\bar f:\overline {\mls  E}\rightarrow \overline {\mls M}_{1,1}$ denote the extension of $\mls  E$ provided by the Tate curve and let $\omega _{\overline {\mls E}/\overline  {\mls M_{1,1}}}$ denote the relative dualizing sheaf.  Then the sheaf $\bar f_*\omega _{\overline {\mls E}/\overline {\mls M_{1,1}}}$ is a line bundle on $\overline {\mls M}_{1,1}$ extending $\lambda $.  In what follows we will abuse notation and write also $\lambda $ for this line bundle on $\overline {\mls M}_{1,1}$.

If $\Lambda $ is a ring and $t:\Sp (\Lambda )\rightarrow \mls M_{1,1}$ is a  morphism corresponding to an elliptic curve $E/\Lambda $, then after replacing $\Lambda $ by an \'etale extension the family $E$ can be described by an equation
\begin{equation}\label{coordinates}
y^2+a_1xy+a_3y=x^3+a_2x^2+a_4x+a_6.
\end{equation}
Define
\begin{equation}
b_2 = a_1^2+4a_2, \ \ b_4 = a_1a_3+2a_4, \ \ b_6 = a_3^2+4a_6, \ \ b_8 = -a_1a_3a_4-a_4^2+a_1^2a_6+a_2a_3^2+4a_2a_6,
\end{equation}
and the discrimant
\begin{equation}\label{discriminant}
\Delta = -b_2^2b_8-8b_4^3-27b_6^2+9b_2b_4b_6\in \Lambda ^*.
\end{equation}
With these chosen coordinates a basis for $t^*\lambda $ is given by the \emph{invariant differential}
\begin{equation}\label{1.0.4}
\pi = dx/(2y+a_1x+a_3).
\end{equation}
Any two choices of coordinates \ref{coordinates} differ by a transformation
\begin{equation}
x' = u^2x+r, \ \ y' = u^3y+su^2x+t,
\end{equation}
where $u\in \Lambda ^*$ and $r, s,t \in \Lambda $.    One can compute that  the invariant differential $\pi '$ obtained from the coordinates $(x', y')$ is equal to $u^{-1}\pi $, and that the discrimant $\Delta '$ in the coordinates $(x', y')$ is equal to $u^{12}\Delta $.  In particular, the element $\Delta \pi ^{\otimes 12}\in t^*\lambda ^{\otimes 12}$ is independent of the choice of coordinates, and therefore defines a  trivialization of $\lambda ^{\otimes 12}$ over $\mls M_{1,1}$.

Let $p :\mls M_{1,1, S}\rightarrow \mathbb{A}^1_S$ be the map defined by the $j$-invariant
\begin{equation}
j=(b_2^2-24b_4)^3/\Delta .
\end{equation}

\begin{thm}\label{thm1} Let $S$  be a  scheme.   Then the map
\begin{equation}\label{1.1.1}
\Z/(12)\times \text{\rm Pic}(\mathbb{A}^1_S)\rightarrow \text{\rm Pic}(\mls M_{1,1,S}), \ \ (i, \mls L)\mapsto \lambda ^{\otimes  i}\otimes p ^*\mls L
\end{equation}
is an isomorphism if either of  the following hold:
\begin{enumerate}
\item [(i)] $S$ is a $\Z[1/2]$-scheme.
\item [(ii)] $S$ is reduced.
\end{enumerate}
\end{thm}

\begin{rem} As we observe in \ref{failremark} the theorem fails for  nonreduced schemes in characteristic $2$.
\end{rem}

\begin{thm}\label{thm2} The map
\begin{equation}
\Z\times  \text{\rm Pic}(S)\rightarrow \text{\rm  Pic}(\overline {\mls M}_{1,1,S}) \ \ (n, M)\mapsto \lambda ^n\otimes _{\mls O_S}M
\end{equation}
is an isomorphism for any scheme $S$.
\end{thm}

\begin{rem} By standard limit arguments it suffices to prove the above results in the case when  $S$ is noetherian.  In what follows we will therefore restrict to the category of noetherian schemes unless otherwise stated.
\end{rem}

\begin{pg}{\bf Acknowledgements.} Olsson partially supported by NSF grant DMS-0555827 and an Alfred P. Sloan fellowship.  
\end{pg}

\section{When $6$ is invertible on $S$}

Though the case when  $6$ is invertible  follows from the more technical work in subsequent sections, we include  here a proof in the case of a $\Z[1/6]$-scheme since it is much easier than the more  general  cases.

Let $\tilde s_4:S\rightarrow \mls M_{1,1,S}$ be the section corresponding to the elliptic curve with automorphism group $\mu _4$ ($y^2=x^3+x$ with $\Delta = -64$, $j=1728$) and $\tilde s_6:S\rightarrow \mls M_{1,1,S}$ the section corresponding to the elliptic curve with automorphism group $\mu _6$ ($y^2+y=x^3$ with $\Delta =-27$, $j=0$).  These sections define closed immersions $s_4:B\mu _{4, S}\hookrightarrow \mls M_{1,1,S}$ and $s_6:B\mu _{6, S}\hookrightarrow \mls M_{1,1,S}$.  For any line bundle $\mls L$ on $\mls M_{1,1,S}$ the pullback $s_4^*\mls L$ (resp. $s_6^*\mls L$) corresponds to a line bundle $M_4$ (resp. $M_6$) on $S$ with action of the group $\mu _4 $ (resp. $\mu _6$).  We thus get maps
$$
\rho _4:\mu _4\rightarrow \underline {\text{Aut}}(M_4) = \mathbb{G}_m, \ \ \rho _6:\mu _6\rightarrow \underline {\text{Aut}}(M_6)\simeq \mathbb{G}_m
$$
defining characters $\chi _4\in \Z/(4)$ and $\chi _6\in \Z/(6)$.

\begin{lem}\label{L1} The pair $(\chi _4, \chi _6)$ lies in $\Z/(12)\subset \Z/(4)\times \Z/(6)$. 
\end{lem}
\begin{proof} The construction of the pair $(\chi _4, \chi _6)$ commutes with arbitrary base change on $S$, so it suffices to consider the case when $S$ is the spectrum of an algebraically closed field, $S = \Sp (k)$.  We have to show that $\rho _4|_{\mu _2} = \rho _6|_{\mu _2}$.  Write $k[[t]]$ for the completion of the local ring of $\mathbb{A}^1_j$ at $j=1728$ and let $k[[z]]$ be the completion of the local ring of $\mls M_{1,1,S}$ at the point coresponding to the curve $y^2=x^3+x$.  Then the map $k[[t]]\rightarrow k[[z]] $ sends $t$ to $z^2$ (with suitable choices of coordinates) and the action of $\mu _4$ is given by $\zeta *z = \zeta ^2\cdot z$.  Write $\mls L|_{k[[z]]]} = k[[z]]\cdot e$ for some basis $e$.  Then $\rho _4$ acts by $\zeta *e = \zeta ^{\chi _4}e$.  From this we see that $\rho _4|_{\mu _2}$ is equal to the character defined by the action of $\mu _2$ on the fiber of $\mls L$ at the generic point of $\mls M_{1,1,S}$.  Similarly, $\rho _6|_{\mu _2}$ is equal to the action on the generic fiber.
\end{proof}

We therefore obtain a map
\begin{equation}\label{retract}
\text{Pic}(\mls M_{1,1,S})\rightarrow \Z/(12), \ \ \mls L\mapsto (\chi _4, \chi _6),
\end{equation}
and it follows from the construction that this map is a homomorphism. Let $K$ denote the kernel.

Recall that a Deligne-Mumford stack $\mls X$ is called \emph{tame} if for every algebraically closed field $\Omega $ and point $\bar x:\Sp (\Omega )\rightarrow \mls X$ the order of the automorphism group of $\bar x$ is relatively prime to the characteristic of $\Omega $.

\begin{lem}\label{L2} Let $\mls X$ be a tame Deligne--Mumford stack with coarse moduli space $\pi :\mls X\rightarrow X$.  Let $\mls L$ be an invertible sheaf on $\mls X$ such that for every geometric point $\bar x\rightarrow \mls X$ the action of the stabilizer group $G_{\bar x}$ on $\mls L(\bar x)$ is trivial.  Then $\pi _*\mls L$ is an invertible sheaf on $X$ and $\pi ^*\pi _*\mls L\rightarrow \mls L$ is an isomorphism.
\end{lem}
\begin{proof} It suffices to prove the lemma after passing to the strict henselization of $X$ at a geometric point $\bar x$.  Let $A = \mls O_{X, \bar x}$ and $B = \mls O_{\mls X, \bar x}$.  Then as explained in \cite[2.12]{Olsson} if $\Gamma $ denotes the stabilizer group of $\bar x$ then there is a natural action of $\Gamma $ on $B$ such that $\mls X = [\Sp (B)/\Gamma ]$.  Let $M$ be the free $B$--module with $\Gamma $--action of rank $1$ defining $\mls L$.  Since $\Gamma $ has order invertible in $k(\bar x)$ (since $\mls X$ is tame) the representation category of $\Gamma $ is semisimple. By our assumptions the reduction $M\otimes k(\bar x)$ is generated by an invariant element and choosing a lifting to an invariant element of $M$ we see that we can write $M = B\cdot e$ where $\Gamma $ acts trivially on $e$.  Then $\pi _*\mls L$ is just $A\cdot e$ and the lemma is immediate.
\end{proof}

\begin{cor} The homomorphism $\pi ^*:\text{\rm Pic}(\mathbb{A}^1_S)\rightarrow K$ is an isomorphism.
\end{cor}
\begin{proof} We show that if $\mls L$ is a line bundle with $(\chi _4, \chi _6) = (0,0)$, then $\pi _*\mls L$ is an invertible sheaf on $\mathbb{A}^1_S$ and $\pi ^*\pi _*\mls L\rightarrow \mls L$ is an isomorphism. By \ref{L2} it suffices to show that for any geometric point $\bar x\rightarrow \mls M_{1,1,S}$ the action of the stabilizer group of $\bar x$ on $\mls L(\bar x)$ is trivial.  For this we may assume that $S$ is the spectrum of an algebraically closed field.  By our assumptions the actions $\rho _4$ and $\rho _6$ are trivial. By the argument used in the proof of \ref{L1} this implies that the action of the generic stabilizer is also trivial.  From this it follows that the action is trivial at every point since over $\mathbb{A}^1-\{0, 1728\}$ the stack $\mls M_{1,1,S}$ is a $\mu _2$--gerbe.
\end{proof}

\begin{lem} The image of $\lambda $ in $\Z/(12)$ is a generator. In particular \ref{retract} is surjective.
\end{lem}
\begin{proof} It suffices to consider the case when $S$ is the spectrum of a field in which case the above shows that $\text{Pic}(\mls M_{1,1,S})$ injects into $\Z/(12)$.  We can in fact compute directly the image of $\lambda $ in $\Z/(4)\times \Z/(6).$ The image in $\Z/(4)$ corresponds to the representation of $\mu _4$ given by the action on the invariant differential $dx/2y$ of the curve $y^2=x^3+x$.  An element $\zeta \in \mu _4$ acts by $(x,y)\mapsto (\zeta ^2x, \zeta y)$ and therefore the action on $dx/2y$ is equal to multiplication by $\zeta $. Therefore the image of $\lambda $ in $\Z/(4)$ is equal to $1$.

Similarly, the image of $\lambda $ in $\Z/(6)$ corresponds to the character given by the invariant differential $dx/(2y+1)$ of the curve $y^2+y=x^3$.  Write $\mu _6 = \mu _2\times \mu _3$.  Then $(-1, 1)$ acts by $(x, y)\mapsto (x, -y-1)$ and $(1, \zeta )$ acts by $(x, y)\mapsto (\zeta x, y)$.  Therefore $(-1, 1)$ acts on the invariant differential by multiplication by $-1$ and $(1,\zeta )$ acts by multiplication by $\zeta $.  It follows that $\lambda $ maps to $1$ in $\Z/(6)$ which implies that $\lambda $ is a generator in $\Z/(12)$.
\end{proof}
\begin{cor} The map $\lambda \times \pi ^*:(\Z/12)\times \text{\rm Pic}(\mathbb{A}^1_S)\rightarrow \text{\rm Pic}(\mls M_{1,1,S})$ is an isomorphism.
\end{cor}

\section{The case of a  normal affine scheme  $S$}\label{section3}

Write $S = \Sp (\Lambda )$  with $\Lambda $ a normal ring.  Let $U$ be the scheme
\begin{equation}
U:= \Sp (\Lambda [a_1, a_2, a_3, a_4, a_6][1/\Delta ]),
\end{equation}
where $\Delta $ is  defined as in \ref{discriminant}.  The equation \ref{coordinates} defines a family of elliptic curves $E\rightarrow U$. Let $G$ denote  the group scheme with underlying scheme $\Sp (\Lambda [u^\pm , r, s, t])$ with group law defined by
\begin{equation}
(u', r', s', t')\cdot (u, r, s, t) = (uu', u^2r'+r, us'+s, u^3t'+u^2r's+t).
\end{equation}
Then $\mls M_{1, 1, S}$ is  isomorphic to the stack theoretic  quotient $[U/G]$.

\begin{prop}\label{normalprop} The  pullback map
\begin{equation}
\text{\rm Pic}(S)\rightarrow \text{\rm Pic}(U)
\end{equation}
is an isomorphism.
\end{prop}
\begin{proof}  The key point is the following result of Ischebeck \cite[\S 4]{Is}.

\begin{lem} Let $\Delta \in \Z[t_1, \dots, t_n]$ be a polynomial satisfying
\begin{enumerate}
\item [(i)] The greatest common divisor of the coefficients of its nonconstant monomials is $1$.
\item [(ii)] For any field $k$ the image of $\Delta $ in $k[t_1, \dots, t_n]$ is irreducible.
\end{enumerate}
Then for any noetherian normal ring $\Lambda $,  the pullback homomorphism
\begin{equation}
\text{\rm  Pic}(\Lambda )\rightarrow \text{\rm Pic}(\Lambda [t_1, \dots, t_n][1/\Delta ])
\end{equation}
is an isomorphism.
\end{lem}
\begin{proof} The assumptions are used as follows:
\begin{enumerate}
\item Assumption (i) implies that the map $\Sp (\Z[t_1, \dots, t_n][1/\Delta ])\rightarrow \Sp (\Z)$ is surjective and hence faithfully flat.  It follows that the map
\begin{equation}
\Sp (\Lambda [t_1, \dots, t_n][1/\Delta ])\rightarrow \Sp  (\Lambda )
\end{equation}
is also faithfully flat.
\item By the preceding observation the divisor $V(\Delta )\subset \Sp (\Lambda [t_1, \dots, t_n])$ does not contain any fibers, and its generic fiber is nonempty and irreducible.  From this it follows that $V(\Delta )$ is irreducible.
\end{enumerate}
It follows that there is an exact sequence of Weil divisor class groups \cite[1.8]{Fulton}
\begin{equation}
\begin{CD}
\Z[V(\Delta )]@>0>> \text{Cl}(\Lambda [t_1, \dots, t_n])@>>> \text{Cl}(\Lambda [t_1, \dots, t_n][1/\Delta ])@>>> 0.
\end{CD}
\end{equation}
We conclude that
\begin{equation}
\text{Cl}(\Lambda )\simeq \text{Cl}(\Lambda [t_1, \dots, t_n])\simeq \text{Cl}(\Lambda [t_1, \dots, t _n][1/\Delta ]).
\end{equation}
The normality of $\Lambda $ implies that the natural  maps from the Picard groups to the Weil divisor class groups are injective.  Thus it suffices to show that if $D\in \text{Cl}(\Lambda )$ is a Weil divisor whose  image in $\text{Cl}(\Lambda [t_1, \dots, t_n][1/\Delta ])$ is in the image of $\text{Pic}(\Lambda [t_1, \dots, t_n][1/\Delta ])$ then $D$ is obtained from a line bundle on $\Sp (\Lambda )$.  This follows from the observation that  $\Lambda \rightarrow \Lambda [t_1, \dots, t_n][1/\Delta ]$ is faithfully flat \cite[\S 4, Satz 6]{Is}.
\end{proof}

We apply the lemma to $\Delta  \in \Z[a_1, \dots, a_6]$.  Then (i) is immediate and (ii) follows from the calculations in \cite[\S 3, \S 4]{Del} (note that though these sections concern characteristics prime to $6$ the same calculations give the irreducibility of $\Delta $ over arbitrary fields).
\end{proof}

The isomorphism $\mls  M_{1,1,S} \simeq [U/G]$ defines a morphism $\sigma :\mls M_{1,1,S}\rightarrow  BG$.  For a character  $\chi :G\rightarrow \mathbb{G}_m$ defining a  line bundle on $BG$ let  $L_\chi  $ be the line  bundle  on $\mls M_{1,1,S}$ obtained by pull back along $\sigma $.

\begin{lem} Let $\mls L$ be a line bundle on $\mls M_{1,1,S}$ such that the pullback $L$ of $\mls L$ to $U$ is trivial.  Then $\mls L\simeq L_\chi $ for some character $\chi :G\rightarrow \mathbb{G}_m$.
\end{lem}
\begin{proof}
Fix a basis $e\in L$.

Let $\mc F$ be the sheaf on the category of affine $S$-schemes (with the \'etale topology) which to any morphism of affine schemes $S'\rightarrow S$ associates $\Gamma (U_{S'}, \mls O_{U_{S'}}^*)$.  There is an inclusion of sheaves $\mathbb{G}_m\subset \mc F$ given by the inclusions $\Gamma (S', \mls O_{S'}^*)\subset \Gamma (U_{S'}, \mls O_{U_{S'}}^*)$.    For any $S'\rightarrow S$ and $g\in G(S')$, we get an element $u_g\in \mc F(S')$ defined by the condition that $g(e) = u_g\cdot e\in L$.  This defines a map of sheaves (not necessarily a homomorphism)
\begin{equation}
f:G\rightarrow \mc F.
\end{equation}
To prove the lemma it suffices to show that $f$ has image contained in $\mathbb{G}_m\subset \mc F$ (note that it is clear that if this holds then the map $G\rightarrow \mathbb{G}_m$ is a homomorphism).

Since $G$ is an affine scheme the map $f$ is determined by a section $u_0\in \mc F(G)$.  Since $G$ is normal and connected, this section $u_0\in \Gamma (U_{G}, \mls O_{U_{G}}^*)$ can be written uniquely as $\beta \Delta ^m$, where $\beta \in \Gamma (G, \mls O_{G}^*)$ and $m\in \Z$.  We need to show that $m=0$.    For this note that the image of  $u_0$ under the map $\mc F(G)\rightarrow \mc F(S)$ defined by the identity section $e:S\rightarrow G$ is equal to $1$.  It follows that $e^*(\beta)\cdot \Delta ^m$ is equal to $1$ in $\Gamma (U, \mls O_U^*)$ which implies that $m=0$. 
\end{proof}

By \ref{normalprop},  if  $\mls L$  is a line bundle on  $\mls M_{1,1,S}$ then the pullback of $\mls L$ to $U$ is  isomorphic to the pullback of  a line bundle $M$ on $S$.  It  follows that any line bundle on $\mls M_{1,1,S}$ is isomorphic to $M\otimes L_{\chi }$ for some character $\chi :G\rightarrow \mathbb{G}_m$.  More such a line bundle $M\otimes L_\chi $ is trivial if and only if $M$ is trivial and $L_\chi $ is trivial.

\begin{lem} Any homomorphism $G\rightarrow  \mathbb{G}_m$ factors through the projection
\begin{equation}\label{3.3.1}
\chi _0:G\rightarrow \mathbb{G}_m, \ \ \ (u, r, s, t)\mapsto u.
\end{equation}
\end{lem}
\begin{proof}
There are three injective homomorphisms
\begin{equation}
j_r, j_s, j_t:\mathbb{G}_a\hookrightarrow G
\end{equation}
sending $x\in \mathbb{G}_a$ to $(1, x, 0,0)$, $(1, 0, x, 0)$, and $(1, 0,0,x)$ respectively.
The formula
\begin{equation}
(1, r, 0,0)(1, 0,s,0)(1, 0,0,t-rs) = (1, r, s, t)
\end{equation}
shows that the subgroup of $G$ generated by the images of these three inclusions is equal to the kernel of $\chi _0$.  
Since any homomorphism $\mathbb{G}_a\rightarrow \mathbb{G}_m$ is trivial, it follows that any homomorphism $G\rightarrow \mathbb{G}_m$ has kernel containing $\text{Ker}(\chi _0)$.
\end{proof}

The line bundle $\lambda $ is  trivialized over $U$ by the invariant differential $\pi $ defined in \ref{1.0.4} and as mentioned in the introduction the action of $(u, r, s, t)\in  G$ on $\pi $ is through the character $G\rightarrow \mathbb{G}_m$ sending $(u,r,s,t)$ to $u^{-1}$.  Putting all this together we find that  \ref{1.1.1} is surjective.

In fact, if $\chi :G\rightarrow \mathbb{G}_m$ is a character, a trivialization of $L_{\chi  }$ is given by a unit $\theta \in \Gamma (U, \mls O_U^*)$ such that for any $(u, r, s, t)\in G$ we have $(u,r,s,t)*\theta  = \chi ^{-1}(u,r,s,t)\theta $.  Any unit $\theta $ on $U$ can be written as $\beta \Delta ^m$ for $\beta \in \Lambda ^*$ and $m\in \Z$.   We have
\begin{equation}
(u,r,s,t)*(\beta \Delta ^m) = \beta u^{12m}\Delta ^m.
\end{equation}
It follows that $L_{\chi }$ is trivial if and only if $\chi  = \chi _0^{12m}$, for some $m$.

This completes the proof of \ref{thm1} in the case  when $S$ is affine and normal. \qed

A very similar argument can be used to prove \ref{thm2} in the case when the base scheme $S$ is affine and normal.  Let $c_4 = b_2^2-24b_4$.  Then one can show (see for example \cite[III.1.4]{Silv}) that \ref{coordinates} is nodal precisely when $\Delta = 0$ and $c_4 \neq 0$.  Let $\widetilde U$ denote 
\begin{equation}
\Sp (\Lambda [a_1, a_2, a_3, a_4, a_6])-V(\Delta, c_4).
\end{equation}
Again the group scheme $G$ acts on $\widetilde U$ and $\overline {\mls M}_{1,1,S}\simeq [\widetilde U/G]$.

\begin{lem} (i) The map 
\begin{equation}
\text{\rm Pic}(\Lambda )\rightarrow \text{\rm Pic}(\widetilde U)
\end{equation}
is an isomorphism.

(ii) The map $\Lambda ^*\rightarrow \Gamma (\widetilde U, \mls O_{\widetilde U}^*)$ is an isomorphism.
\end{lem}
\begin{proof} Statement (ii) is immediate.  Statement (i) follows from a very similar argument to the proof of \ref{normalprop}.  The only new ingredient is that the polynomial $c_4$ is not irreducible over fields of characteristics $2$ and $3$ but it is a power of an irreducible polynomial (in characteristic $2$ it is equal to $a_1^4$ and in characteristic $3$ is is equal to $(a_1+a_2)^2$.
\end{proof}

Using this one sees as before that the map
\begin{equation}
\text{Pic}(\Lambda )\times \text{Pic}(BG)\rightarrow \text{Pic}(\overline {\mls M}_{1,1,S})
\end{equation}
is an isomorphism with the character \ref{3.3.1} mapping to $\lambda ^{-1}$.

\section{The case when $S$ is reduced}

\begin{pg} If $S$ is an arbitrary scheme, and $\mls L$ a line  bundle on $\mls M_{1,1,S}$ then there is a unique function $s\mapsto l(s)\in \Z/(12)$ which associates to a point $s$ the unique power $l(s)$ of $\lambda $ such that $\mls L_s\otimes \lambda ^{-l(s)}$ on $\mls M_{1,1,k(s)}$  descends to $\mathbb{A}^1_{k(s)}$.
\end{pg}

\begin{lem} The function $s\mapsto  l(s)$ is a locally constant function on $S$.
\end{lem}
\begin{proof}
The assertion is local on $S$ so we may assume that $S$ is affine. Furthermore, the assertion can be verified on  each irreducible component so we may assume that $S$ is integral.  Finally if $\widetilde S\rightarrow S$ is the normalization then it suffices to verify the assertion for $\widetilde S$.  In this case the result follows from section \ref{section3}.
\end{proof}

\begin{pg} In particular if $S$ is connected we obtain a homomorphism
\begin{equation}\label{degmap}
\text{Pic}(\mls M_{1,1,S})\rightarrow \Z/(12)
\end{equation}
sending $\lambda $ to $1$.
Thus in general to prove \ref{thm1} we need to show that the kernel of \ref{degmap} is isomorphic to $\text{Pic}(\mathbb{A}^1_S)$.
\end{pg}

\begin{lem} For any locally noetherian scheme $S$, the map $\pi :\mls M_{1,1,S}\rightarrow \mathbb{A}^1_S$ given by the $j$-invariant identifies $\mathbb{A}^1_S$ with the coarse moduli space of $\mls M_{1,1,S}$.
\end{lem}
\begin{proof}
Let $\tilde \pi :\mls M_{1,1,S}\rightarrow X$ be the coarse moduli space (which exists by \cite{KM}).  By the universal property of the coarse moduli space, there exists a unique morphism $f:X\rightarrow \mathbb{A}^1_j$ such that $f\circ \tilde \pi = \pi $.  Since $\pi $ is proper and quasi-finite, the morphism $f$ is also proper and quasi-finite and therefore $f$ is finite.  Furthermore, by \cite{KM} we have  $\tilde \pi _*\mls O_{\mls M_{1,1,S}} = \mls O_{X}$.  It therefore suffices to show that the map $\mls O_{\mathbb{A}^1_j}\rightarrow \pi _*\mls O_{\mls M_{1,1,S}}$ is an isomorphism. It suffices to verify this locally in the flat topology on $S$, so we may further assume that $S$ is the spectrum of a complete noetherian local ring $A$.  In addition, since the morphism $\pi $ is proper, the theorem on formal functions for stacks \cite[\S 3]{Chow} implies that it suffices to show the result over $\Sp (A/\mathfrak{m}_A^n)$ for all $n$.  This reduces the proof to the case when $S$ is the spectrum of an artinian local ring $A$.  Let $k$ be the residue field of $A$, and let $J\subset A$ be an ideal with $J$ annihilated by the maximal ideal of $A$ (so that $J$ is a $k$-vector space).  Set $A_0:= A/J$.  Pushing forward the exact sequence
\begin{equation}
0\rightarrow J\otimes \mls O_{\mls M_{1,1,k}}\rightarrow \mls O_{\mls M_{1,1,A}}\rightarrow \mls O_{\mls M_{1,1,A_0}}\rightarrow 0
\end{equation}
to $\mathbb{A}^1_A$ we obtain a commutative diagram
\begin{equation}
\begin{CD}
0@>>> (\pi _*\mls O_{\mls M_{1,1,k}})\otimes J@>>> \pi _*\mls O_{\mls M_{1,1,A}}@>>> \pi _*\mls O_{\mls M_{1,1,A_0}}@. \\
@. @AaAA @AbAA @AAcA @. \\
0@>>> (\mls O_{\mathbb{A}^1_k})\otimes J@>>> \mls O_{\mathbb{A}^1_A}@>>> \mls O_{\mathbb{A}^1_{A_0}}@>>> 0.
\end{CD}
\end{equation}
By induction and the case when $A$ is a field, we get that $a$ and $c$ are isomorphisms and therefore $b$ is an isomorphism also.
\end{proof}

\begin{pg} To complete the proof of \ref{thm1} in the case when $S$ is reduced, we make some general observations about the relationship between line bundles on a stack and line bundles on the coarse moduli space.

Let $S$ be a  noetherian scheme and $\mls X\rightarrow S$ a Deligne-Mumford stack over $S$.  Let $\pi :\mls X\rightarrow X$ be the coarse moduli space, and assume that the formation of the coarse space $X$ commutes with arbitrary base change on $S$ and that $X$ is reduced (we just saw that this holds for $\mls M_{1,1}$ over a reduced scheme).  For a field valued point $x:\Sp (k)\rightarrow S$ let $\pi _x:\mls X_x\rightarrow X_x$ denote the base change $\mls X\times _Sx\rightarrow X\times _Sx$.
\end{pg}

\begin{prop} Let $L$ be a line bundle on $\mls X$ such that for every field valued point $x:\Sp (k)\rightarrow S$ the sheaf $\pi _{x*}(L|_{{\mls X}_x})$ is locally free of rank $1$ and  $\pi _x^*\pi _{x*}(L|_{{\mls X}_x})\rightarrow L|_{{\mls X}_x}$ is an isomorphism.  If $\mls X\rightarrow X$ is flat, then the sheaf $\pi _*L$ is locally free of rank $1$ on $X$ and $\pi ^*\pi _*L\rightarrow L$ is an isomorphism.
\end{prop}
\begin{proof}
One immediately reduces to the case when $X = \Sp (R)$, $Y = \Sp (B)$ is a finite flat $R$ scheme, $\Gamma $ is a finite group acting on $Y$ over $X$ such that $\mls X = [Y/\Gamma ]$ (indeed \'etale locally on the coarse space every Deligne-Mumford stack can be presented in this way \cite[2.12]{Olsson}).  Let $M$ denote the $B$--module corresponding to $L$, so that $M$ comes equipped with an action of $\Gamma $ over the action on $B$. We can even assume that $R$ is a local ring and  that  $M$ is a free $R$--module (forgetting the $B$--module structure). 
We are then trying to compute the kernel of the map
$$
M\rightarrow \prod _{\gamma \in \Gamma }M, \ \ m\mapsto (\cdots, \gamma (m)-m, \cdots )_{\gamma \in \Gamma }.
$$
We can also assume that $S = \Sp (\Lambda )$ is affine.
\end{proof}

\begin{lem} Let $R$ be a reduced local $\Lambda $--algebra   and let $A\in M_{n\times m}(R)$  be a matrix (which we view as a map $R^n\rightarrow R^m$)  with the property that for every $x\in \Sp (\Lambda )$ the matrix $A(x)\in M_{n\times m}(R\otimes _\Lambda k(x))$ has kernel a free $R\otimes _\Lambda k(x)$--space of rank $1$.  Then $\text{\rm Ker}(A)$ is a free rank $1$ module over $R$ and for every $x\in \Sp (\Lambda )$ the natural map $\text{\rm Ker}(A)\otimes _\Lambda k(x)\rightarrow \text{\rm Ker}(A(x))$ is an isomorphism.
\end{lem}
\begin{proof} By induction on $n$.  If $n=1$, then the assertion is that $A$ is a matrix with $A(x)$ the zero matrix for all $x\in \Sp (\Lambda )$.  Since $R$ is reduced this implies that $A$ is the zero matrix.

For the inductive step consider the system of $m$ equations
$$
\sum _{i}a_{ij}X_i = 0
$$
that we are trying to solve in $R$.  If $x\in \Sp (\Lambda )$ is the image of the closed point of $\Sp (R)$, then $A(x)$ is not zero since $n\geq 2$.  Since $R$ is local some $a_{ij}$ is invertible and so we can solve for the variable $X_i$.  This gives a system of $m-1$--equations in $n-1$ variables which again has the property that for every point $x\in \Sp (\Lambda )$ the image in $R\otimes k(x)$ has a unique line of solutions.  By induction we obtain the result.   
\end{proof}

This completes the proof of (\ref{thm1} (i)).

\section{Proof of (\ref{thm1} (ii))}\label{section5}

\begin{prop}\label{5.1} For any scheme $S$ over $\Z[1/2]$  and any coherent $\mls O_S$--module $M$, the sheaf $R^1\pi _*(\mls O_{\mls M_{1,1, S}}\otimes _{\mls O_S}M)$ is zero, where  $\pi :\mls M_{1,1, S}\rightarrow \mathbb{A}^1_{j, S}$ is the projection.
\end{prop}
\begin{proof} Using the theorem of formal functions one is reduced to the case when $S$ is the spectrum of a field. Furthermore, if the characteristic is not $3$ the result is immediate, so it suffices to consider $S = \Sp (k)$ with $\text{char}(k) = 3$, and $M = k$.  We may further assume that $k$ is algebraically closed.

The   coherent sheaf $R^1\pi _*(\mls O_{\mls M_{1,1, k}})$ restrict to the zero sheaf on $\mathbb{A}^1_k-\{0\}$, since over this open subset of $\mathbb{A}^1_k$ the stack $\mls M_{1,1,k}$ is tame (the automorphism groups are $\{\pm 1\}$).  Let $\bar x\rightarrow \mls M_{1,1,k}$ be a geometric point mapping to $0$ in $\mathbb{A}^1_k$, and let $A$ denote the completion of  $\mls O_{\mls M_{1,1,k}, \bar x}$ along the maximal ideal. Let $\Gamma _{\bar x}$ denote the stabilizer group scheme of $\bar x$, so that $\Gamma _{\bar x}$ acts on $A$.  The ring of invariants $B:=A^{\Gamma _{\bar x}}$ is equal to the completion  of $\mathbb{A}^1_{k}$ at the origin.  Let  $F$ denote the finite type $B$-module obtained by pulling back $R^1\pi _*(\mls O_{\mls M_{1,1, k}})$ to $\Sp (B)$.  Then $F$ is equal to the cohomology group $H^1(\Gamma _{\bar x}, A)$.  We show that this group is zero.  
Since $F$ is supported on the closed point of $\Sp (B)$, there exists an integer $n$ such that $j^nF = 0$ (where $j\in B$ is the uniformizer defined by the standard coordinate on $\mathbb{A}^1$).  To prove the proposition it therefore suffices to show that $F$ is $j$-torsion free.

For this we use an explicit description of $A$ and $\Gamma _{\bar x}$ given by the Legendre family.  Let
\begin{equation}
V = \Sp (k[\lambda ][1/\lambda (\lambda -1)])
\end{equation}
and let $E_V\rightarrow V$ be the elliptic curve
\begin{equation}
E_V: \ Y^2Z = X(X-Z)(X-\lambda Z).
\end{equation}
If $\mu $ denotes $\lambda +1$, then the $j$-invariant of $E_V$ is equal to $\mu ^6/(\mu ^2-1)^2$ (recall that $\text{char}(k)=3$).  The map $V\rightarrow \mls M_{1,1,k}$ defined by $E_V$ is \'etale, so this defines an isomorphism $A\simeq k[[\mu ]]$.  The group $\Gamma _{\bar x}$ sits in an exact sequence
\begin{equation}
1\rightarrow \{\pm 1\}\rightarrow \Gamma _{\bar x}\rightarrow S_3\rightarrow 1,
\end{equation}
and the action of $\Gamma _{\bar x}$ on $A\simeq k[[\mu]]$ factors through the action of $S_3$ on $k[[\mu ]]$ given by the two automorphisms
\begin{equation}
\alpha : \ \mu \mapsto -\mu, 
\end{equation}
and
\begin{equation}
\beta : \ \mu \mapsto \mu /(1-\mu) = \mu (1+\mu+\mu^2+\cdots ).
\end{equation}
Also note that the Leray spectral sequence
\begin{equation}
E_2^{pq} = H^p(S_3, H^q(\{\pm 1\}, A))\implies H^{p+q}(\Gamma _{\bar x}, A)
\end{equation}
and the fact that $H^q(\{\pm 1\}, A) = 0$ for $q>0$ (since $2$ is invertible in $k$) implies that $H^1(\Gamma _{\bar x}, A) = H^1(S_3, A)$.

An element in $H^1(S_3, A)$ can be represented by a set map $\xi :S_3\rightarrow k[[\mu ]]$ (written $\sigma \mapsto \xi _\sigma $) such that for $\sigma, \tau \in S_3$ we have (recall the action is a right action)
\begin{equation}\label{cohcondition}
\xi _{\sigma \tau } = \xi _\sigma ^\tau +\xi _\tau .
\end{equation}
The class of $\xi $ is trivial if there exists an element $g\in k[[\mu ]]$ such that $\xi _\sigma = g^\sigma -g$ for all $\sigma \in S_3$.  Note that \ref{cohcondition} implies that it suffices to check the equalities $\xi _\sigma = g^\sigma -g$ for a set of generators $\sigma \in S_3$.

If $\xi $ represents a class in $H^1(S_3, A)$ annihilated by $j$, there exists an element $g\in k[[\mu ]]$ such that 
\begin{equation}
\frac{\mu ^6}{(\mu ^2-1)^2}\xi _\sigma = g^\sigma -g
\end{equation}
for all $\sigma \in S_3$.  To prove that $H^1(S_3, A)$ is $j$-torsion free, it therefore suffices to show that for such a $\xi $ we can choose $g$ to have $\mu$-adic valuation $\geq 6$ (since $A$ is $j$-torsion free).

For this note that we can without loss of generality assume that $g$ has no constant term, and then write
\begin{equation}
g = a_1\mu + a_2\mu ^2+a_3\mu ^3+a_4\mu ^4+a_5\mu ^5+g_{\geq 6},
\end{equation}
where $g_{\geq 6}$ has $\mu$-adic valuation $\geq 6$.  We have
\begin{equation}
\frac{\mu ^6}{(\mu ^2-1)^2}\xi _\alpha = 2a_1\mu +2a_3\mu ^3+2a_5\mu ^5+(g_{\geq 6}^\alpha -g_{\geq 6})
\end{equation}
which implies that $a_1 = a_3 = a_5 = 0$.  Then
\begin{equation}
\frac{\mu ^6}{(\mu ^2-1)^2}\xi _\beta = 2a_2\mu ^3+(\text{higher order terms})
\end{equation}
which gives $a_2=0$.  Finally using this we see that
\begin{equation}
\frac{\mu ^6}{(\mu ^2-1)^2}\xi _\beta = a_4\mu ^5+(\text{higher order terms})
\end{equation}
which implies that $a_4 = 0$ as desired.
This completes the proof of \ref{5.1}.
\end{proof}

\begin{pg} Now let us prove \ref{thm1} for a connected $\Z[1/2]$-scheme $S$.
We need to show  that if $L$ is a line bundle on $\mls M_{1,1,S}$ such that for any field valued point  $s\in S$ the fiber $L_s$ on $\mls M_{1, 1, s}$ descends to $\mathbb{A}^1_{j, s}$ then $L$ descends to $\mathbb{A}^1_{j, S}$.  By a standard limit argument it suffices to consider the case when $S$ is noetherian and even affine, say $S = \Sp (\Lambda )$.  Let $J\subset \Lambda $ denote the nilradical.  We would like to  inductively show that if the result holds for over $\Lambda /J^r$ then it also holds for $\Lambda /J^{r+1}$.  In other words, let $L_0$ denote  a line bundle on $\mathbb{A}^1_{j, \Lambda /J^r}$ and $\widetilde L$ a lifting of $\pi ^*L_0$ to $\mls M_{1, 1, \Lambda /J^{r+1}}$.  Then we want to show that $\widetilde L$ is pulled back from a lifting of $L_0$ to $\mathbb{A}^1_{j, \Lambda /J^{r+1}}$.  By standard deformation theory this is equivalent to showing that the map
$$
0 = H^1(\mathbb{A}^1_{\Lambda }, J^r/J^{r+1}) \rightarrow H^1(\mls M_{1,1,\Lambda }, J^r/J^{r+1})
$$
is an isomorphism.  Equivalently that $H^1(\mls M_{1,1, \Lambda }, J^r/J^{r+1})$ is zero.   Since $\mathbb{A}^1_{j, \Lambda }$ is affine, the group $H^1(\mls M_{1,1, \Lambda }, J^r/J^{r+1})$ is zero if and only if the sheaf $R^1\pi _*(J^r/J^{r+1}\otimes \mls O_{\mls M_{1,1, \Lambda }})$ is zero on $\mathbb{A}^1_{j, \Lambda }$ which follows from \ref{5.1}.
This completes the proof of \ref{thm1}. \qed
\end{pg}

\section{Computations  in characteristic $2$}

\begin{prop} Let $k$ be a field of characteristic $2$, and let $\pi :\overline {\mls M}_{1,1,k}\rightarrow \mathbb{P}^1_k$ be the morphism defined by the $j$-invariant.  Then $R^1\pi _*\mls O_{\overline {\mls M}_{1,1,k}}$ is a line bundle on $\mathbb{P}^1_k$ of negative degree.
\end{prop}
\begin{proof}
We may without loss of generality assume that $k$  is algebraically closed.

Let $\mls U_\infty \subset \overline {\mls M}_{1,1,k}$ denote the open substack $\overline {\mls M}_{1,1,k}\times _{\mathbb{P}^1_j}\mathbb{A}^1_{1/j}$ (the complement of $j=0$), and let $\mls U_0 = \mls M _{1,1,k}\subset \overline {\mls M}_{1,1,k}$ denote the complement of $j = \infty $.  Let $U_\infty , U_0\subset \mathbb{P}^1_j$ be the coarse moduli spaces (the standard open cover of $\mathbb{P}^1_j$).

The stack $\mls U_\infty $ is a $\Z/(2)$--gerbe over $U_\infty $.  Now in general, if $f:\mls G\rightarrow X$ is a $\Z/(2)$--gerbe in characteristic $2$, the sheaf $R^1f_*\mls O_{\mls G}$ is locally free of rank $1$ and in fact canonically trivialized.  This can be seen as follows.  Etale locally on $X$, we have $\mls G = X\times B(\Z/(2))$.  Computing in this local situation, one sees that $R^1f_*(\Z/(2))$ is a locally constant sheaf of groups \'etale locally isomorphic to $\Z/(2)$, and the natural map $R^1f_*(\Z/(2))\otimes _{\Z/(2)}\mls O_X\rightarrow R^1f_*\mls O_{\mls G}$ (which exists since we are in characteristic $2$) is an isomorphism.  Since a group of order $2$ admits no nontrivial automorphisms there is a canonical isomorphism $\Z/(2)\simeq R^1f_*(\Z/(2))$ which induces a canonical trivialization of $R^1f_*\mls O_{\mls G}$.  In the case of $\mls G = X\times B(\Z/(2))$ and $X = \Sp (A)$ we have
$$
H^1(\mls G, \mls O_{\mls G}) \simeq \text{Hom}_{\text{Gp}}(\Z/(2), A)
$$
and the trivialization is given by the homomorphism sending $1\in \Z/(2)$ to $1\in A$.

\begin{lem} The sheaf $R^1\pi _*\mls O_{\overline {\mls M}_{1,1,k}}$ is locally free of rank $1$ on $\mathbb{P}^1_j$.
\end{lem}
\begin{proof}
By finiteness of coherent cohomology for stacks the sheaf is in any case coherent.  Since $\mathbb{P}^1_j$ is a smooth curve it therefore suffices to show that it is torsion free.  Furthermore, the only issue is at the point $j=0$.  Since the formation of cohomology commutes with flat base change, it suffices to show that 
\begin{equation}\label{cohgroup}
H^1(\overline {\mls M}_{1,1,k}\times _{\mathbb{P}^1_j}\Sp (k[[j]]), \mls O_{\overline {\mls M}_{1,1,k}\times _{\mathbb{P}^1_j}\Sp (k[[j]])})
\end{equation}
is $j$--torsion free. 

For this we use the so-called Hesse presentation of $\mls M_{1,1,k}$.  Let
\begin{equation}
V = \Sp (k[\mu, \omega][1/(\mu ^3-1)]/(\omega ^2+\omega +1)),
\end{equation}
and let $E_V\rightarrow V$ be the elliptic curve given by the equation 
\begin{equation}
X^3+Y^3+Z^3 = \mu XYZ.
\end{equation}
This is elliptic curve has a basis for its three-torsion group given by the points $[1:0:-1]$ and $[-1:\omega :0]$.  In fact, this is the universal elliptic curve with full level three structure.  The $j$-invariant of $E_V$ is $\mu ^{12}/(\mu ^3-1)^3$ (since we are in characteristic $2$).  In particular, the fiber  over $j=0$ is the curve $X^3+Y^3+Z^3 = 0$.

Changing the choice of basis for the $3$-torsion subgroup defined an action of $GL_2(\mathbb{F}_3)$ on $V$ such that $\mls M_{1,1,k}\simeq [V/GL_2(\mathbb{F}_3)].$  A calculation shows that this action is described as follows:
\begin{enumerate}
\item $(\mu, \omega )*\begin{pmatrix} 1 & 0 \\ -1 & 1\end{pmatrix} = (\omega \mu, \omega ).$
\item $(\mu, \omega )*\begin{pmatrix} 0 & -1 \\ 1 & 0 \end{pmatrix} = (\mu /(\mu -1), \omega ).$
\item $(\mu, \omega )*\begin{pmatrix} 1 & 0 \\ 0 & -1 \end{pmatrix} = (\mu, \omega ^2).$
\end{enumerate}

Putting this together one finds that
\begin{equation}\label{Hesse}
\overline {\mls M}_{1,1,k}\times _{\mathbb{P}^1_j}\Sp (k[[j]])\simeq [\Sp (k[[\mu ]])/\text{SL}_2(\mathbb{F}_3)],
\end{equation}
where $\alpha = \begin{pmatrix} 1 & 0 \\ -1 & 1\end{pmatrix}$ acts by $\mu \mapsto \zeta \mu $ (for some fixed primitive cube root of unity $\zeta $) and $\beta = \begin{pmatrix} 0 & -1 \\ 1 & 0 \end{pmatrix}$ acts by $\mu \mapsto \mu /(\mu -1)$.  

As in the proof of \ref{5.1}, an element of \ref{cohgroup} is given by a set map $\xi :\text{SL}_2(\mathbb{F}_3)\rightarrow k[[\mu ]]$ (written $\sigma \mapsto \xi _\sigma $) such that for any two elements $\sigma , \tau \in \text{SL}_2(\mathbb{F}_3)$ we have 
$$
\xi _{\sigma \tau } = \xi _{\sigma }^\tau +\xi _\tau ,
$$
and the class of $\xi $ is trivial if there exists an element $g\in k[[\mu ]]$ such that for every $\sigma $ we have $\xi _\sigma = g^\sigma -g$.

Now if \ref{cohgroup} has $j$--torsion there exists a set map $\xi $ as above and an element $g\in k[[\mu ]]$ such that for all $\sigma $ we have
$$
\frac{\mu ^{12}}{\mu ^3-1}\xi _\sigma  = g^\sigma -g.
$$
To prove that \ref{cohgroup} is torsion free it suffices to show that we can choose $g$ to be divisible by $\mu ^{12}$.  For since $k[[\mu ]]$ is an integral domain we then have
$$
\xi _\sigma = (\frac{\mu ^3-1}{\mu ^{12}}g)^\sigma - (\frac{\mu ^3-1}{\mu ^{12}}g).
$$

We can without loss of generality assume that $g$ has no constant term.  Write
$$
g = a_1\mu + a_2\mu ^2+\cdots + a_{11}\mu ^{11} + g_{\geq 12}.
$$
Then $g^{\alpha } - g$ has $\mu $-adic valuation $\geq 12$ (recall that $\alpha = \begin{pmatrix} 1 & 0 \\ -1 & 1\end{pmatrix}$).  Looking at the coefficients $a_i$ this implies that all but $a_3$, $a_6$, and $a_9$ are zero, so 
$$
g = a_3\mu ^3+a_6\mu ^6+a_9\mu ^9+g_{\geq 12}.
$$
Similarly $g^\beta - g$ has $\mu $--adic valuation $\geq 12$.  Looking at the coefficient of $\mu ^4$ in $g^\beta -g$ one sees that $a_3 = 0$.  Then looking at the coefficent of $\mu ^8$ one sees that $a_6=0$, and finally looking at the coefficient of $\mu ^{10}$ one sees that $a_9=0$.
\end{proof}

 Let $M$ denote the cohomology group \ref{cohgroup} (a $k[[j]]$--module) and let $M_\eta $ denote $M\otimes _{k[[j]]}k[[j]][1/j]$.  Let $e_\infty \in M_\eta $ denote the basis element defined by the canonical trivialization of $R^1\pi _*\mls O_{\overline {\mls M}_{1,1,k}}$ over $\mls U_\infty $.  The lattice $M\subset M_\eta $ defines a valuation $\nu $ on $M_\eta $ and it suffices to show that $\nu (e_\infty )<0$.  Equivalently we have to show that for any element $m\in M$ if we write $m = he_\infty $ in $M_\eta $ then the $j$--adic valuation of $h$ is positive.

For this we again use the presentation \ref{Hesse}.  An element $m\in M$ is then represented by a map $\xi :\text{SL}_2(\mathbb{F}_3)\rightarrow k[[\mu ]]$.  The corresponding element in $M_\eta $ can in terms of the basis $e_\infty $ be described as follows.  First of all the element $\xi _{\beta ^2}\in k[[\mu ]]$ is then $\text{SL}_2(\mathbb{F}_3)$--invariant since for any other element $\sigma $ we have
$$
\xi _{\beta ^2}^\sigma + \xi _\sigma = \xi _{\beta ^2\sigma } = \xi _{\sigma \beta ^2} = \xi _\sigma ^{\beta ^2}+\xi _{\beta ^2}
$$
and $\beta ^2$ acts trivially on $k[[\mu ]]$.  Therefore $\xi _{\beta ^2}$ is actually an element in $k[[j]]$.  The image of $\xi $ in $M_\eta \simeq \text{Hom}(\Z/(2), k[[j]][1/j])$ is then equal to the homomorphism 
$$
\Z/(2)\rightarrow k[[j]][1/j], \ \ 1\mapsto \xi _{\beta ^2}.
$$
The class $e_\infty $ corresponds to the homomorphism sending $1$ to $1$ so we have to show that the $j$--adic valuation of $\xi _{\beta ^2}$ is positive.  For this let $f = \xi _\beta $.  Then 
$$
\xi _{\beta ^2} = f^\beta +f = f(\mu (1+\mu +\mu ^2+\cdots ))+f(\mu ).
$$
Since we are in characteristic $2$ it follows that the $\mu $--adic valuation of $\xi _{\beta ^2}$ is at least $2$, and therefore the $j$--adic valuation of $\xi _{\beta ^2}$ is also positive.
\end{proof}

\begin{cor} For any field $k$, we have $H^1(\overline {\mls M}_{1,1,k}, \mls O_{\overline {\mls M}_{1,1,k}}) = 0$.
\end{cor}
\begin{proof} We have $R^1\pi _*\mls O_{\overline {\mls M}_{1,1,k}}  = 0$ when $\text{char}(k) \neq 2 $ (when $\text{char}(k)=3$ this follows from section \ref{section5}).  It follows that
\begin{equation}
H^0(\mathbb{P}^1_k, R^1\pi _*\mls O_{\overline {\mls M}_{1,1,k}}) = 0
\end{equation}
in all characteristics.
 From the Leray spectral sequence we obtain
\begin{equation}
H^1(\overline {\mls M}_{1,1,k}, \mls O_{\overline {\mls M}_{1,1,k}}) = H^1(\mathbb{P}^1_k, \mls O_{\mathbb{P}^1_k}) = 0.
\end{equation}
\end{proof}

\begin{rem}\label{failremark} Note that if $\text{char}(k) = 2$, then the restriction of $R^1\pi _*\mls O_{\mls M_{1,1,k}}$ to $\mathbb{A}^1_k\subset \mathbb{P}^1_k$ is non-zero.  From the Leray spectral sequence it follows that the map
\begin{equation}
0 = H^1(\mathbb{A}^1_k, \mls O_{\mathbb{A}^1_k})\rightarrow H^1(\mls M_{1,1,k}, \mls O_{\mls M_{1,1,k}})
\end{equation}
is \emph{not} an isomorphism.  Since the group $H^1(\mls M_{1,1,k}, \mls O_{\mls M_{1,1,k}})$ classifies deformations of the structure sheaf to $\mls M_{1,1,k[\epsilon ]/(\epsilon ^2)}$ this implies that there are line bundles on $\mls M_{1,1,k[\epsilon ]/(\epsilon ^2)}$ which are in  the kernel of \ref{degmap} but are nontrivial.  This implies that \ref{thm1} fails for $S = \Sp (k[\epsilon ]/(\epsilon ^2))$. More generally, \ref{thm1} fails for any nonreduced affine scheme over $\mathbb{F}_2$.
\end{rem}

\section{Proof  of  \ref{thm2}}

\begin{pg} In order to prove \ref{thm2} it is easiest to prove a stronger statement that implies it.  Let  $\mls Pic(\overline {\mls M}_{1,1,S})$ denote the Picard stack over $S$ which to any $S$--scheme $T$ associates the group if line bundles on $\overline {\mls M}_{1,1,T}$. By \cite[5.1]{Aoki},  the stack $\mls Pic(\overline {\mls M}_{1,1,S})$ is an algebraic stack (an Artin stack) over $S$.   There is a morphism of stacks
\begin{equation}\label{picstackmap}
\Z\times B\mathbb{G}_{m, S}\rightarrow \mls Pic(\overline {\mls M}_{1,1,S})
\end{equation}
sending a pair $(n, L)$ consisting of $n\in \Z$ and $L$ a line bundle on  $S$ to $\lambda ^n\otimes _{\mls O_S}L$ on $\overline {\mls M}_{1,1,S}$.  The following theorem implies \ref{thm2} by passing evaluating both sides of \ref{picstackmap} on $S$ and passing to isomorphism classes.
\end{pg}

\begin{thm}\label{thm3}
The morphism \ref{picstackmap} is an isomorphism.
\end{thm}
\begin{proof}
Note first that if $n$ and $n'$ are integers and $L$ and $L'$ are line bundles on  $S$, then  $\lambda ^n\otimes L$ and $\lambda ^{n'}\otimes L'$ on $\overline {\mls M}_{1,1,S}$ are isomorphic if and only if $n=n'$ and $L\simeq L'$.  Indeed, if these two sheaves are isomorphic, then this implies that $\lambda ^{n-n'}$ descends to $\mathbb{P}^1_j$.  By the case of a field this implies that $n=n'$.  In this case we recover $L$ and $L'$ from their pullbacks to $\overline {\mls M}_{1,1,S}$ by pushing back down to $S$.  Therefore,
the functor \ref{picstackmap} is  fully faithful. It therefore suffices to show that for any cartesian diagram
\begin{equation}
\begin{CD}
P@>>> S\\
@VVV @VVLV \\
\Z\times B\mathbb{G}_m@>>> \mls Pic _{\overline {\mls M}_{1,1,S}}
\end{CD}
\end{equation}
the morphism of algebraic spaces $P\rightarrow S$ is an isomorphism. For this it suffices to consider the case when $S$ is artinian local.  Furthermore, we know the result in the case when $S$ is the spectrum of a field by section \ref{section3}. Since a line bundle on the spectrum of an artinian local ring is trivial, what we therefore need to show is that if $S$ is an artinian local ring then any line bundle on $\overline {\mls M}_{1,1,S}$ is isomorphic to $\lambda ^n$ for some $n$.  Proceeding by induction on the length of $S$, it further suffices to consider the following.  Let $S = \Sp (A)$, $k$ the residue field of $A$, and let $J\subset A$ be a square--zero ideal annihilated by the maximal ideal of $A$, and set $A_0 = A/J$.  Then any deformation of $\lambda ^n$ over $\overline {\mls M}_{1,1, A_0}$ to $\overline {\mls M}_{1,1,A}$ is isomorphic to $\lambda ^n$.  Using the exponential sequence
$$
0\rightarrow J\otimes \mls O_{\overline {\mls M}_{1,1,k}}\rightarrow \mls O^*_{\overline {\mls M}_{1,1,A}}\rightarrow \mls O_{\overline {\mls M}_{1,1,A_0}}^*\rightarrow 0
$$
one sees that this amounts exactly to $H^1(\overline {\mls M}_{1,1,k}, \mls O_{\overline {\mls M}_{1,1,k}}) = 0$.
\end{proof}

\providecommand{\bysame}{\leavevmode\hbox
to3em{\hrulefill}\thinspace}

\end{document}